\documentclass[11pt]{amsart}
\usepackage[T1]{fontenc}
\usepackage{amsmath,amssymb,amsthm}
\usepackage{etoolbox}
\makeatletter
\patchcmd{\@setauthors}{\MakeUppercase{\authors}}{\authors}{}{%
  \PackageError{authorcase}{Could not patch amsart author capitalization}{}%
}
\makeatother
\usepackage[margin=1.1in]{geometry}
\usepackage[hidelinks]{hyperref}
\hypersetup{
  pdftitle={A Higher-Dimensional Extension of the Laugesen--Siudeja Lower Bound},
  pdfauthor={Zhengchao Ji}
}

\numberwithin{equation}{section}
\newtheorem{theorem}{Theorem}[section]
\newtheorem{lemma}[theorem]{Lemma}
\newtheorem{proposition}[theorem]{Proposition}
\newtheorem*{conjecture}{Conjecture}
\newtheorem*{problem}{Problem}
\theoremstyle{remark}
\newtheorem*{remark}{Remark}

\newcommand{\R}{\mathbb R}
\newcommand{\dd}{\,\mathrm d}
\DeclareMathOperator{\Cov}{Cov}
\DeclareMathOperator{\tr}{tr}
\DeclareMathOperator{\Var}{Var}

\title[A higher-dimensional Laugesen--Siudeja bound]
{A Higher-Dimensional Extension of the Laugesen--Siudeja Lower Bound}
\author[Zhengchao Ji]{%
Zhengchao Ji%
\texorpdfstring{%
\protect\\[5pt]
\protect\normalfont\protect\small
Department of Mathematics, China Jiliang University%
\protect\\
\protect\normalfont\protect\small
Hangzhou 310018, People's Republic of China%
\protect\\
\protect\normalfont\protect\small
Correspondence: jizhengchao@zju.edu.cn%
}{}%
}
\date{}

\begin{document}

\begin{abstract}
We prove a higher-dimensional analogue of the lower bound in the
Laugesen--Siudeja conjecture, using a covariance normalization of
the first Dirichlet eigenvalue. We give a covariance-dependent
lower bound for every bounded convex domain. From this bound we obtain
the optimal constant for convex planar domains and for simplices in
every dimension \(d\geq2\). If the simplex covariance inequality
stated in a 2026 preprint holds, the same constant is optimal over all
bounded convex domains in dimensions \(d\geq3\). In each of these cases,
the infimum is approached by compressing a regular simplex in the
direction of an edge and is unattained.
\end{abstract}

\maketitle

\section{Introduction and main result}

For \(d\geq2\), we write \(\mathcal C_d\) for the class of nonempty
bounded convex open subsets of \(\R^d\). Given \(\Omega\in\mathcal C_d\),
we consider the Dirichlet eigenvalue problem
\[
 -\Delta u=\lambda u\quad\text{in }\Omega,
 \qquad u=0\quad\text{on }\partial\Omega
\]
An eigenvalue is a value of \(\lambda\) for which this problem has
a nonzero weak solution \(u\in H^1_0(\Omega)\). We write
\(0<\lambda_1(\Omega)<\lambda_2(\Omega)\leq\cdots\) for the
eigenvalues, counted with multiplicity. The first satisfies
\[
 \lambda_1(\Omega)
 =\inf_{0\neq u\in H^1_0(\Omega)}
 \frac{\int_\Omega |\nabla u|^2\dd x}{\int_\Omega |u|^2\dd x}.
\]

We write \(V=|\Omega|\) for the volume and
\[
 c=\frac1V\int_\Omega x\dd x,\qquad
 I=\int_\Omega |x-c|^2\dd x
\]
for the centroid and the polar moment of inertia about the centroid.
In dimension two we write \(A=V\) for the area and use the
scale-invariant functional
\[
 Q(\Omega)=\lambda_1(\Omega)\frac{A^3}{I}.
\]
Indeed, a dilation by \(r\) multiplies \(\lambda_1\), \(A\), and \(I\) by
\(r^{-2}\), \(r^2\), and \(r^4\), respectively.

P\'olya and Szeg\H{o}'s inclusion method gives upper and lower bounds
for \(Q\) on convex domains~\cite[Sections~1.19 and~5.11b]{PSzego}.
Laugesen and Siudeja used rotational symmetry and tight frames to
show that \((\lambda_1+\cdots+\lambda_n)A^3/I\) is maximized by the
equilateral triangle among triangles and by the square among
parallelograms~\cite{LS}. For \(n=1\), the equilateral triangle and
every rectangle give the same value \(12\pi^2\). This observation
led to the following conjecture.

\begin{conjecture}[Laugesen--Siudeja, {\cite[Conjecture~4.1]{LS}}]
For every nonempty bounded convex planar domain \(\Omega\),
\[
 \frac92\pi^2<Q(\Omega)\leq12\pi^2.
\]
Equilateral triangles and rectangles give equality on the right.
Thin acute isosceles triangles and thin circular sectors approach
the lower endpoint.
\end{conjecture}

Freitas obtained more precise eigenvalue asymptotics for thin
triangles~\cite{FreitasThin}. The conjecture was restated in
Laugesen's survey~\cite[Conjecture~7.1]{LaugesenSurvey} and in the
chapter by Laugesen and Siudeja on special
domains~\cite[Conjecture~6.62]{LSChapter}.

In our earlier work~\cite[Theorems~1.1 and~1.2]{JiPlanar}, we settled
the planar conjecture. The lower bound holds for every bounded convex
planar domain and is sharp along a family of thin acute isosceles
triangles. For the upper bound, we used Temple's inequality and exact
arithmetic to show that an explicit isosceles trapezoid has
\(Q>12\pi^2\). We now ask for the corresponding sharp lower bound in higher
dimensions. To formulate it, we first rewrite \(Q\) in covariance
notation.

For \(\Omega\in\mathcal C_d\), we write \(C\) for the covariance
matrix of a point uniformly distributed in \(\Omega\). Thus
\[
 C=\frac1V\int_\Omega (x-c)(x-c)^{\mathsf T}\dd x.
\]
The matrix \(C\) is positive definite and satisfies \(I=V\tr C\).
In dimension two, the identity \(\tr(C^{-1})=\tr C/\det C\) expresses
\(Q\) in terms of the trace of \(C^{-1}\) and the determinant of \(C\).
We use the following covariance form of the normalization in dimension \(d\)
\begin{equation}\label{eq:normalization}
 \widetilde Q_d(\Omega)
 =\frac{\lambda_1(\Omega)}{\tr(C^{-1})}
   \left(\frac{V^2}{\det C}\right)^{2/d},
\end{equation}
which is invariant under similarities and agrees with \(Q\) in
dimension two.

This is the first-eigenvalue normalization of Laugesen and
Siudeja~\cite[Corollary~3.3]{LSHigher} written in terms of the
covariance of \(\Omega\). Indeed, write \(\Omega=T(D)\) after
translation, where \(T\) is invertible and \(D\) has covariance
\(\sigma^2\mathrm{Id}_d\) for some \(\sigma>0\).
Since \(C_\Omega=\sigma^2TT^{\mathsf T}\) and
\(C_{T^{-\mathsf T}D}=\sigma^2T^{-\mathsf T}T^{-1}\),
the volume factors and the identity \(I=V\tr C\) give
\[
 \widetilde Q_d(TD)
 =\lambda_1(TD)\,|TD|^{2/d}
   \frac{|T^{-\mathsf T}D|^{1+2/d}}{I(T^{-\mathsf T}D)}.
\]
We study sharp lower bounds for this functional.

For a rectangular box with side lengths \(a_1,\ldots,a_d\), we have
\[
 \lambda_1=\pi^2\sum_{i=1}^d a_i^{-2},\qquad
 C=\frac{1}{12}\operatorname{diag}(a_1^2,\ldots,a_d^2),
 \qquad V^2=\prod_{i=1}^d a_i^2.
\]
As in the box example in \cite[Section~3]{LSHigher}, these formulas
give \(\widetilde Q_d=12\pi^2\) for every such box.
For the scale-invariant quantity \(\lambda_1 V^{1+4/d}/I\), a box
with side lengths \(L,1,\ldots,1\) gives
\[
 12\pi^2\,\frac{L^{4/d}(L^{-2}+d-1)}{L^2+d-1},
\]
which tends to zero as \(L\to\infty\) for \(d>2\).

For the lower endpoint of \(\widetilde Q_d\), we set
\begin{equation}\label{eq:sharp-constant}
 c_d=\frac{\pi^2(d+2)(d+1)^{1+2/d}}{2(d!)^{4/d}}.
\end{equation}

Motivated by the thin triangles in the plane, we compress a regular
simplex in the direction of an edge. Along this family,
\(\widetilde Q_d\) tends to \(c_d\). This leads to the following
problem.

\begin{problem}[Sharp lower bound in higher dimensions]
For \(d\geq3\), determine
\[
 \inf_{\Omega\in\mathcal C_d}\widetilde Q_d(\Omega).
\]
Is this infimum equal to \(c_d\) in \eqref{eq:sharp-constant}?
Can any \(\Omega\in\mathcal C_d\) attain it?
\end{problem}

In the plane, we obtain the covariance estimate from the published
random-triangle inequality of Campi, Colesanti, and
Gronchi~\cite[\S3, inequality~(11)]{CCG}. For general convex domains
in dimensions \(d\geq3\), we use the simplex covariance inequality
stated in the preprint by OpenAI~\cite[Theorem~1.1 and
equation~(7.6)]{SimplexPreprint}. We have not independently verified
the proof of this inequality. We therefore present the
higher-dimensional argument separately from \cite{JiPlanar} and
assume the covariance inequality only for the sharp constant over
all convex domains in dimensions \(d\geq3\).

\begin{theorem}\label{thm:main}
Let \(d\geq2\).
\begin{enumerate}
\item For every \(\Omega\in\mathcal C_d\),
\begin{equation}\label{eq:det-bound}
 \widetilde Q_d(\Omega)>
 \frac{\pi^2}{2(d+1)(d+2)}
 \left(\frac{V^2}{\det C}\right)^{2/d}.
\end{equation}
The coefficient in \eqref{eq:det-bound} is optimal. For every
open \(d\)-simplex \(S\), \eqref{eq:det-bound} gives
\(\widetilde Q_d(S)>c_d\).
The infimum over open \(d\)-simplices is \(c_d\), approached by
compressing a regular simplex in the direction of one of its edges.

\item If \(d=2\), or if \(d\geq3\) and the simplex covariance
inequality \eqref{eq:simplex-covariance} holds for every
\(d\)-dimensional convex body, then every \(\Omega\in\mathcal C_d\)
satisfies
\begin{equation}\label{eq:sharp-lower}
 \widetilde Q_d(\Omega)>c_d.
\end{equation}
In either case,
\[
 \inf_{\Omega\in\mathcal C_d}\widetilde Q_d(\Omega)=c_d.
\]
The same family of compressed simplices approaches this infimum,
which is unattained.
\end{enumerate}
\end{theorem}

\begin{remark}
The hypothesis \eqref{eq:simplex-covariance} is used only to pass
from the covariance-dependent bound \eqref{eq:det-bound} to the
constant \(c_d\) over all convex domains in dimensions \(d\geq3\).
It is the inequality stated in \cite[equation~(7.6)]{SimplexPreprint}.
If we take that result as established, the hypothesis may be omitted
from the second conclusion of Theorem~\ref{thm:main}.
\end{remark}

Formula \eqref{eq:sharp-constant} gives \(c_2=9\pi^2/2\) and
\(c_3=10\pi^2/3^{4/3}\).

In Section~\ref{sec:variance} we establish the sharp variance bound
for one-dimensional projections. We apply this bound in the principal
directions of the covariance matrix in Section~\ref{sec:spectral},
where it yields a strict spectral estimate. Finally, in
Section~\ref{sec:sharpness} we compute the covariance determinant of
a simplex and complete the proof of Theorem~\ref{thm:main} by
constructing a family that approaches the lower bound.

\section{A sharp variance bound for projected measures}\label{sec:variance}

\begin{lemma}\label{lem:variance}
Let \(p>0\), and suppose that the support of a probability density
\(f\) is an interval of length \(W\). If \(f^{1/p}\) is concave on the
interior of that interval, then
\begin{equation}\label{eq:variance}
 \Var(f)\geq\frac{W^2}{2(p+2)(p+3)}.
\end{equation}
The bound is attained by the density proportional to
\((1-2|t|/W)^p\) on \([-W/2,W/2]\).
\end{lemma}

\begin{proof}
We first reduce to a symmetric density. Write \(m\) for the mean of
\(f\) and \(f^*\) for its symmetric decreasing rearrangement on
\([-a,a]\), where \(a=W/2\). Among sets of length
\(\ell_u=|\{f>u\}|\), the interval centered at \(m\) minimizes the
integral of \((x-m)^2\). Applying this comparison at each level \(u\),
we obtain
\begin{align}\label{eq:rearrangement}
 \Var(f)
 &=\int_0^\infty\int_{\{f>u\}}(x-m)^2\dd x\dd u \notag\\
 &\geq \int_0^\infty\int_{-\ell_u/2}^{\ell_u/2}x^2\dd x\dd u
 =\int_{-a}^{a} x^2 f^*(x)\dd x .
\end{align}

To check that rearrangement preserves concavity of the \(p\)-th root,
we put \(g=f^{1/p}\) and \(I_s=\{g>s\}\) for \(0\leq s<\sup g\).
If \(r=(1-\theta)s+\theta t\) with \(0<\theta<1\), then
\[
 (1-\theta)I_s+\theta I_t\subseteq I_r
 \quad\Longrightarrow\quad
 (1-\theta)|I_s|+\theta|I_t|\leq |I_r|.
\]
The superlevel sets \(J_s=\{g^*>s\}\) of \(g^*=(f^*)^{1/p}\)
are centered intervals with \(|J_s|=|I_s|\). Thus the length
inequality gives \((1-\theta)J_s+\theta J_t\subseteq J_r\).
For \(x,y\in(-a,a)\), we let \(s\to g^*(x)\) and
\(t\to g^*(y)\) from below to obtain concavity of \(g^*\). Since its
one-sided limit at \(a\) is nonnegative, we have
\[
 \frac{g^*(x)}{a-x}\leq\frac{g^*(y)}{a-y},
 \qquad 0\leq x<y<a.
\]
We now compare \(f^*\) with the density
\[
 f_0(x)=\frac{p+1}{2a}(1-|x|/a)^p,\qquad |x|<a.
\]
The preceding inequality shows that \(f^*/f_0\) is nondecreasing
on \([0,a)\). Since both densities have mass \(1/2\) on this interval,
we can choose \(x_0\in[0,a]\) so that
\((x^2-x_0^2)(f^*-f_0)(x)\geq0\) almost everywhere on \([0,a]\).
Using symmetry and equality of the masses, we get
\[
 \int_{-a}^{a} x^2(f^*-f_0)(x)\dd x
 =2\int_0^a(x^2-x_0^2)(f^*-f_0)(x)\dd x\geq0.
\]
Substituting in \eqref{eq:rearrangement}, we obtain
\[
 \Var(f)\geq\Var(f_0)
 =(p+1)a^2\int_0^1 t^2(1-t)^p\dd t
 =\frac{2a^2}{(p+2)(p+3)}
 =\frac{W^2}{2(p+2)(p+3)}.
\]
The density \(f_0\) attains the bound.
\end{proof}

\section{The spectral comparison}\label{sec:spectral}

\begin{proposition}\label{prop:spectral}
Every \(\Omega\in\mathcal C_d\) satisfies
\begin{equation}\label{eq:spectral}
 \lambda_1(\Omega)>
 \frac{\pi^2}{2(d+1)(d+2)}\tr(C^{-1}).
\end{equation}
Consequently, \eqref{eq:det-bound} holds.
\end{proposition}

\begin{proof}
We first choose an orthonormal eigenbasis \(e_1,\ldots,e_d\) with
\(Ce_i=s_i e_i\). Write \(W_i\) for the width in direction \(e_i\)
and \(R\) for the open supporting box in these directions.
To estimate \(W_i\), we use the density of the projection onto \(e_i\),
\[
 f_i(t)=\frac{1}{V}
 \mathcal H^{d-1}\bigl(\Omega\cap\{x:x\cdot e_i=t\}\bigr).
\]
Here \(\mathcal H^{d-1}\) measures the \((d-1)\)-dimensional volume
of the section. By Brunn--Minkowski in dimension \(d-1\),
\(f_i^{1/(d-1)}\) is concave on the interior of its support.
We can therefore apply Lemma~\ref{lem:variance} to get
\[
 s_i=\Var(f_i)\geq\frac{W_i^2}{2(d+1)(d+2)}.
\]
We now compare \(\Omega\) with \(R\). Domain monotonicity and the
box spectrum give
\[
 \lambda_1(\Omega)\geq\lambda_1(R)
 =\pi^2\sum_{i=1}^d W_i^{-2}
 \geq \frac{\pi^2}{2(d+1)(d+2)}\tr(C^{-1}).
\]

It remains to justify strictness. If \(\Omega\subsetneq R\), then
\(R\setminus\overline\Omega\) contains an open ball, since
\(\Omega=\operatorname{int}\overline\Omega\).
Equality of the first eigenvalues would make the zero extension of
a first eigenfunction on \(\Omega\) a minimizer of the Rayleigh
quotient on \(R\). Every first eigenfunction of \(R\) is a nonzero
multiple of a product of first sine modes and is nowhere zero in
\(R\). The ball on which the extension vanishes thus forces
\(\lambda_1(\Omega)>\lambda_1(R)\).
When \(\Omega=R\), we instead use the box formulas directly
\[
 \frac{\lambda_1(R)}{\tr(C^{-1})}
 =\frac{\pi^2}{12}>\frac{\pi^2}{2(d+1)(d+2)}.
\]
We obtain \eqref{eq:spectral} in both cases. Dividing by
\(\tr(C^{-1})\) and multiplying by \((V^2/\det C)^{2/d}\)
gives \eqref{eq:det-bound}.
\end{proof}

\section{The simplex constant and its sharpness}\label{sec:sharpness}

\begin{lemma}\label{lem:det}
Every \(d\)-simplex \(S\) satisfies
\begin{equation}\label{eq:det-simplex}
 \frac{\det\Cov(X_S)}{|S|^2}
 =\frac{(d!)^2}{(d+1)^{d+1}(d+2)^d},
\end{equation}
where \(X_S\) is uniform on \(S\).
\end{lemma}

\begin{proof}
We write \(C_S=\Cov(X_S)\) and first record how the determinant
changes under an affine map. If \(S'=BS+b\) with \(B\) invertible, then
\[
 C_{S'}=BC_SB^{\mathsf T},
 \qquad |S'|=|\det B|\,|S|,
 \qquad \frac{\det C_{S'}}{|S'|^2}=\frac{\det C_S}{|S|^2}.
\]
We may therefore take \(S=\operatorname{conv}(0,e_1,\ldots,e_d)\),
where the \(e_i\) are the coordinate unit vectors. By integrating
the coordinate monomials, we obtain
\[
 |S|=\frac1{d!},\qquad
 \mathbb E X_i=\frac1{d+1},\qquad
 \mathbb E(X_iX_j)=\frac{1+\delta_{ij}}{(d+1)(d+2)}.
\]
With \(\boldsymbol 1=(1,\ldots,1)^{\mathsf T}\), it follows that
\[
 C_S=\frac{1}{(d+1)(d+2)}
 \left(\mathrm{Id}_d-\frac{\boldsymbol 1\boldsymbol 1^{\mathsf T}}{d+1}\right).
\]
The matrix in parentheses acts by \(1/(d+1)\) on
\(\operatorname{span}\{\boldsymbol 1\}\) and by the identity on its
orthogonal complement. Hence
\[
 \frac{\det C_S}{|S|^2}
 =\frac{(d!)^2}{[(d+1)(d+2)]^d(d+1)}.
\]
\end{proof}

\begin{proof}[Proof of Theorem~\ref{thm:main}]
\emph{Lower bounds.}
Proposition~\ref{prop:spectral} gives \eqref{eq:det-bound} in every
dimension. We now derive the constant bound \(c_d\) from
\eqref{eq:det-bound} by estimating the covariance determinant. We take \(K=\overline\Omega\),
so \(|K|=V\) and \(C_K=C\) because the boundary has measure zero.
For \(d\geq3\), part~(2) assumes the covariance bound stated
in \cite[equation~(7.6)]{SimplexPreprint},
\begin{equation}\label{eq:simplex-covariance}
 \frac{\det C}{V^2}
 \leq \frac{(d!)^2}{(d+1)^{d+1}(d+2)^d}.
\end{equation}
For \(d=2\), we recover this bound from random triangles.
Take independent uniform points \(X_0,X_1,X_2\) in \(K\), set
\(Y_i=X_i-c\), and write \([x,y]=\det[x,y]\). Then twice the signed
area of the resulting triangle is
\[
 D=[Y_0,Y_1]+[Y_1,Y_2]+[Y_2,Y_0].
\]
By independence, \(\mathbb E\bigl([Y_i,Y_j]^2\bigr)=2\det C\) for
\(i\ne j\). Conditioning on the shared point in each mixed product,
we find that its expectation is zero. Thus, with
\(T=\operatorname{conv}\{X_0,X_1,X_2\}\),
\[
 \frac{\mathbb E|T|^2}{V^2}
 =\frac{\mathbb E D^2}{4V^2}
 =\frac{6\det C}{4V^2}
 =\frac{3\det C}{2V^2}.
\]
The normalized moment \(\mathbb E|T|^2/V^2\) is unchanged when
\(K\) is scaled to unit area. Applying the planar random-triangle
inequality \cite[\S3, inequality~(11)]{CCG} with \(r=2\) and \(n=3\)
shows that the normalized moment is bounded above by its value on a
triangular domain. Lemma~\ref{lem:det} then gives
\eqref{eq:simplex-covariance} for \(d=2\).

Under the hypotheses of part~(2), combining \eqref{eq:det-bound}
and \eqref{eq:simplex-covariance} gives
\begin{align*}
 \widetilde Q_d(\Omega)
 &> \frac{\pi^2}{2(d+1)(d+2)}
     \left(\frac{V^2}{\det C}\right)^{2/d}\\
 &\geq\frac{\pi^2}{2(d+1)(d+2)}
     \left(\frac{(d+1)^{d+1}(d+2)^d}{(d!)^2}\right)^{2/d}
 =c_d.
\end{align*}
For a simplex, we obtain the same strict bound by substituting
\eqref{eq:det-simplex} directly into \eqref{eq:det-bound}.

\emph{Sharpness and nonattainment.}
We now construct a family approaching the lower bound. We translate
a regular open simplex \(S\) so that an edge of \(\overline S\)
has endpoints \(v_\pm=\pm We/2\), where \(e\) is a unit vector and
\(W\) is the edge length. The other vertices lie in \(e^\perp\).
Writing \(S_t=\overline S\cap\{x:x\cdot e=t\}\), we have, for
\(0\leq t<W/2\) and \(x\in S_t\) with barycentric coefficients
\(\alpha_\pm\) at \(v_\pm\),
\[
 \alpha_+-\alpha_-=\frac{2t}{W},
 \qquad
 S_t=\frac{2t}{W}v_++\left(1-\frac{2t}{W}\right)S_0.
\]
To obtain the section identity, we subtract \(2t/W\) from the
coefficient of \(v_+\) and renormalize the remaining coefficients.
Convexity gives the reverse inclusion. By reflection, we therefore have
\[
 \mathcal H^{d-1}(S_t)
 =\left(1-\frac{2|t|}{W}\right)^{d-1}\mathcal H^{d-1}(S_0),
 \qquad |t|<W/2.
\]
Reflection across \(e^\perp\) preserves \(S\), so \(e\) is a
covariance eigenvector. The normalized section volumes give the
density \(f_0\) from Lemma~\ref{lem:variance} with \(p=d-1\).
Hence we obtain
\begin{equation}\label{eq:edge-var}
 C_Se=se,\qquad s=\frac{W^2}{2(d+1)(d+2)}.
\end{equation}

Next we compress \(S\) by a factor \(\varepsilon\) along \(e\) and
write \(S_\varepsilon\) for the resulting simplex. To compare
eigenvalues, we fix \(0<\delta<1\) and \(z\in S\). The segment with
endpoints \((1-\delta)v_\pm+\delta z\) is compactly
contained in \(S\). Hence \(I_\delta\times B_\delta\subset S\)
in coordinates along and perpendicular to \(e\), where \(I_\delta\)
has length \((1-\delta)W\) and \(B_\delta\) is a transverse box.
We also enclose \(S\) in a box of width \(W\) in direction \(e\).
After compressing both boxes, we use domain monotonicity to obtain
\[
 \frac{\pi^2}{W^2}
 \leq \varepsilon^2\lambda_1(S_\varepsilon)
 \leq \frac{\pi^2}{(1-\delta)^2W^2}
      +\varepsilon^2\lambda_1(B_\delta).
\]
We first let \(\varepsilon\to0^+\) with \(\delta\) fixed and
then let \(\delta\to0^+\). It follows that
\begin{equation}\label{eq:thin-spectral}
 \lim_{\varepsilon\to0^+}\varepsilon^2\lambda_1(S_\varepsilon)
 =\frac{\pi^2}{W^2}.
\end{equation}

Compression replaces the eigenvalue \(s\) in direction \(e\) by
\(\varepsilon^2s\) and preserves the perpendicular eigenvalues.
The ratio \(|S|^2/\det C_S\) is affine invariant by
Lemma~\ref{lem:det}. Using \eqref{eq:edge-var} and
\eqref{eq:thin-spectral}, we get
\begin{align*}
 \varepsilon^2\tr(C_{S_\varepsilon}^{-1}) &\longrightarrow s^{-1},
 \qquad
 \frac{|S_\varepsilon|^2}{\det C_{S_\varepsilon}}
 =\frac{|S|^2}{\det C_S},\\
 \widetilde Q_d(S_\varepsilon)
 &\longrightarrow \frac{\pi^2s}{W^2}
     \left(\frac{|S|^2}{\det C_S}\right)^{2/d}
 =c_d.
\end{align*}
The limit in \eqref{eq:thin-spectral}, together with
\(\varepsilon^2\tr(C_{S_\varepsilon}^{-1})\to s^{-1}\), shows that
the coefficient of \(\tr(C^{-1})\) in \eqref{eq:spectral} is
optimal. Equivalently, the coefficient in \eqref{eq:det-bound} is
optimal without the covariance hypothesis. The limit of
\(\widetilde Q_d(S_\varepsilon)\) gives \(c_d\) as the infimum over
simplices and, under the hypotheses of the second conclusion, over
\(\mathcal C_d\). The strict bounds give nonattainment.
\end{proof}

\section*{Declaration of competing interests}
The author declares no competing interests.

\section*{Data availability}
No datasets were used in this study.

\section*{Funding}
This work was supported by the National Natural Science Foundation of
China (grant number 12501068) and by the Natural Science Foundation of
Zhejiang Province (grant number LQN25A010001).

\section*{Acknowledgements}
The author is grateful to Professors Kefeng Liu and Hongwei Xu for their
continuous support, valuable advice, and encouragement.

\noindent\textit{AI assistance.}
The author used OpenAI Codex (GPT-6.1 sol) to explore the higher-dimensional
covariance normalization, develop and check parts of the mathematical
argument, and assist with English-language editing and LaTeX
typesetting. The author takes responsibility for the final manuscript.

\end{document}